\documentclass{article}
\usepackage[utf8]{inputenc}
\usepackage[a4paper,top=3cm,bottom=2cm,left=2cm,right=2cm]{geometry}
\usepackage{amsmath}
\usepackage{amssymb} 
\usepackage{amsthm}
\usepackage{mathtools}
\usepackage{enumerate}
\usepackage{comment}
\usepackage{graphicx}
\usepackage{braket}
\usepackage{booktabs}
\usepackage{url}
\usepackage[affil-it]{authblk}

\usepackage{listings}
\usepackage{xcolor}
\definecolor{codegreen}{rgb}{0,0.6,0}
\definecolor{codegray}{rgb}{0.5,0.5,0.5}
\definecolor{codepurple}{rgb}{0.58,0,0.82}
\definecolor{backcolour}{rgb}{0.98,0.98,0.95}

\lstdefinestyle{mystyle}{
  backgroundcolor=\color{backcolour},   commentstyle=\color{codegreen},
  keywordstyle=\color{magenta},
  numberstyle=\tiny\color{codegray},
  stringstyle=\color{codepurple},
  basicstyle=\ttfamily\footnotesize,
  breakatwhitespace=false,         
  breaklines=true,                 
  captionpos=b,                    
  keepspaces=true,                 
  numbers=left,                    
  numbersep=5pt,                  
  showspaces=false,                
  showstringspaces=false,
  showtabs=false,                  
  tabsize=2
}

\usepackage[square,numbers]{natbib}
\newtheorem*{phm}{Phantom}

\title{Exact sinogram: an analytical approach\\ to the Radon transform of phantoms}
\author{Monica Dessole  \textsuperscript{1}, Marta Gatto  \textsuperscript{1}, Davide Poggiali  \textsuperscript{2}, Francesca Tedeschi  \textsuperscript{1}}

\affil{ \footnotesize
	  \textsuperscript{1} Department of Mathematics ``Tullio Levi Civita'', University of Padova, Italy. \\
	  \textsuperscript{2} PNC - Padova Neuroscience Center, Italy.   \medskip}
\date{}

\begin{document}

\maketitle

\begin{abstract}
Phantoms can serve as a gold standard for the validation of MRI numerical methods. In some special cases, it is possible to compute analytically the Radon transform, or sinogram, of a phantom. In this work, we present analytical formulae to compute the exact sinograms of three classes of phantoms. We compare the use of the discrete Radon transform, that yields an approximate sinogram, and the correspondent analytical sinogram for image reconstruction.  
\end{abstract}

\section{Introduction and motivations}
The imaging technique of exposing an object with rays at different angles is called \emph{tomography}, and it produces a projection image of the inaccessible regions of the body. The regions are determined by their \textit{attenuation function}, that quantifies how the intensity of radiations is attenuated, proportionally to the density, due to mass absorption.
The inverse problem of recovering a function from its integral values along cross sections was first studied by Johann Radon in 1917~\cite{4307775}, implemented by \citet{Cormack} to Computerized Tomography (CT) technologies in the 1960s.
The graphic representation of the Radon transform is commonly referred to as \emph{sinogram}. \\
Numerical methods for image reconstruction are tested on the so-called \emph{phantoms}, fictitious images that reproduce a body section, as the Shepp-Logan phantom \cite{shepplogan}. Using a phantom gives the advantage to test a reconstruction algorithm on a zero-noise data so the error we get is only due to numerical inaccuracies in the algorithm itself~\cite{Feeman, KaipioSomersalo, Hansen}.  
In some special cases, it is possible to compute analytically the Radon transform of a phantom. \\
The scope of this work is to compare the use of the discrete Radon transform, that yields an approximate sinogram, and the correspondent analytical sinogram for image reconstruction. \\
We created \texttt{exact-sinogram} \cite{exact_sinogram}, an open-source python library of Analytical Radon Transforms for different classes of phantoms, obtained as a combination of ellipses, squares or rectangles, in order to successively reconstruct the original image and to analyze the error.

\section{Mathematical tools}
Let $f = f(x,y)$ an integrable continuous function with compact support, standing for the attenuation function, with $(x,y) \in \Omega \subset \mathbb{R}^2$. We define the \textit{Radon transform} of $f$, as the following functional
\begin{equation*}
    Rf(t,\theta):= \int_{l_{t,\theta}} f \,ds = \int_{s=-\infty}^{\infty} f(t\cos{\theta}-s\sin{\theta},t\sin{\theta}+s\cos{\theta}) \,ds. \qquad \qquad \forall \, (t,\theta) \in \mathbb{R} \times [0,2\pi),
\end{equation*}
where $l_{t,\theta}$ is the line that passes through the point $(t\cos(\theta), t \sin(\theta))$ and orthogonal to to the unit vector $\mathbf{n} = (\cos(\theta), \sin(\theta))$.
The classical way of representing the Radon transform is by treating $(t,\theta)$ as rectangular coordinates. The values of the Radon transform are depicted in a continuum grey scale from $0$ to $1$, standing for black and white respectively. This picture is the sinogram and it represents data generated by the X-ray emission/detection machine for the given slice of the sample.

We aim to reconstruct an approximation of the function $f$ from a set of values of the Radon transform calculated on a discrete grid $(t, \theta)$, using mathematical phantoms.
A phantom is an image in grey scale with colors representing mass absorption, made of a collection of geometric figures (possibly convex), each of which has a different constant attenuation. 

In this paper we exploit the analytical Radon transform of two types of figures: the ones that can be described by a single mathematical equation, e.g. ellipses, and the polygons, e.g. squares and rectangles.

We give here a brief description of two examples. More details can be found in \cite{Feeman}.

\begin{phm}[Phantom of ellipses]
Every ellipse in the $xy$-plane is uniquely determined by six real parameters: its semi-axes $a$ and $b$, the coordinates ($x_0, y_0$) of its center point, the angle $\phi$ of rotation with respect to the $x$-axis, and the value $\delta$ of the attenuation function in the region of the plane inside the ellipse. The coordinates in the roto-translated frame are $\hat{x}=(x-x_0)\cos{\phi}+(y-y_0)\sin{\phi}$ and $\hat{y}=(y-y_0)\cos{\phi}-(x-x_0)\sin{\phi}$.\\
Now, let $\epsilon$ be a generic ellipse and for real numbers $t$ and $\theta$, let pose $\hat{\theta}=\theta-\phi$ and $\hat{t}=t-x_0\cos{\theta}-y_0\sin{\theta}$. 
We assume the attenuation function has the following form
\begin{equation*}
    f_{\epsilon}(x,y)\coloneqq 
    \begin{cases}
    \delta & \text{if } \frac{\hat{x}^2}{a^2} + \frac{\hat{y}^2}{b^2} \leq 1,\\
    0      &  \text{otherwise}.
    \end{cases}
\end{equation*}
The Radon transform of $f_{\epsilon}$ is obtained computing analytically the intersections between the line $l_{t,\theta}$ and the ellipse, making the difference between them, i.e. their distance, and is expressed as 
\begin{equation}
\label{eq:radon_ellipse}
\begin{cases}
    Rf_{\epsilon}(t,\theta)= \delta \frac{2ab\sqrt{b^2\sin^2{\hat{\theta}}+a^2\cos^2{\hat{\theta}}-\hat{t}^2}}{b^2\sin^2{\hat{\theta}}+a^2\cos^2{\hat{\theta}}} & \text{if }  \hat{t}^2\leq b^2\sin^2{\hat{\theta}}+a^2\cos^2{\hat{\theta}},\\
    Rf_{\epsilon}(t,\theta)=0 & \text{otherwise}.
    \end{cases}
\end{equation}
\end{phm}

\begin{phm}[Phantom of squares]
A square in general position in the plane is, instead, defined by five parameters: coordinates $(x_0,y_0)$ of the center, side length $w$, angle of counterclockwise rotation $\phi$ from the horizontal and the attenuation value $\delta$.  
To generate the sinogram, we imagine the region inside the square as the intersection of two infinite bands: one vertical and one horizontal. Then, we compute the points at which the line crosses these bands to find, among them, the ``entry" and ``exit" of the line within the square. The distance between the exit and the entry points will be the Radon value associated with those $(t,\theta)$, analogously to the case of the ellipse.
\end{phm}
When two planar figures overlap in a phantom, the attenuation coefficient on the overlap will be the sum of the individual attenuation coefficients. The linearity property of the Radon transform will allow us to compute the analytical sinogram of any combination of such figures. 

For each phantom, we are going to calculate two reconstructed images, one from the approximated sinogram, computed with the discrete Radon transform, and another from the exact sinogram. To do this, we apply the inverse Radon transform provided by the iradon function of the Python library Scikit-image \cite{scikitimage}. We would expect a smaller error with the analytical approach. 

\begin{figure}
    \centering
    \includegraphics[width=0.32\textwidth, trim=1cm 0cm 1cm 0cm]{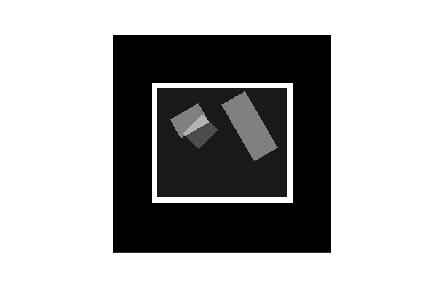} 
    \includegraphics[width=0.57\textwidth, trim=1.2cm 7.8cm 1cm 7.8cm,clip]{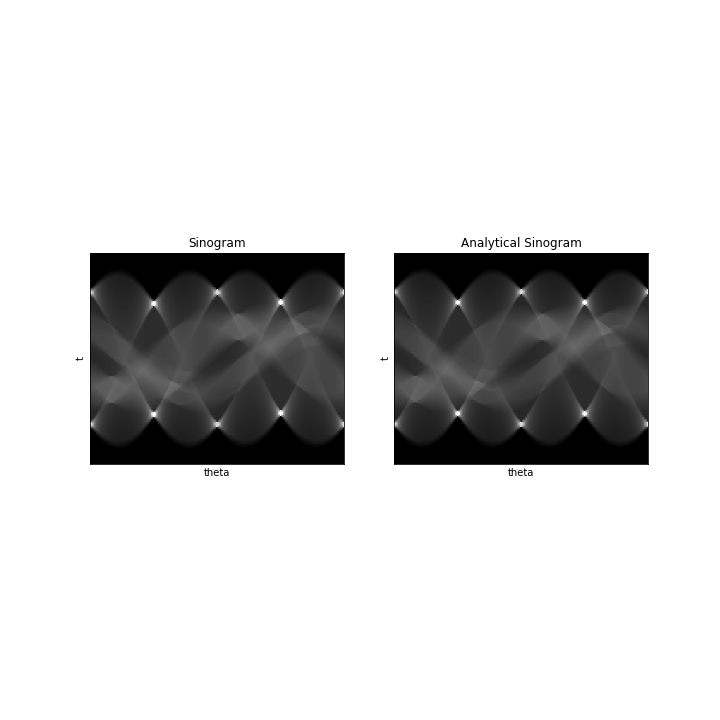}
    \caption{Phantom of rectangles (left), its approximated (center) and exact (right) sinograms.}
    \label{fig:phm_ellipse}
\end{figure}

\section{The library}
The python library \texttt{exact-sinogram} is available in Github (\url{https://github.com/francescat93/Exact_sinogram}) or can be installed via pip. It is written on top of numpy~\cite{van2011numpy}, scipy~\cite{2020SciPy-NMeth}, scikit-image~\cite{scikitimage}, and matplotlib for displaying the images.\\

A phantom is an instance of the class \verb|Phantom|, initialized with the phantom type (\verb|`ellipses'|, \verb|`squares'|, \verb|`rectangles'|) and a matrix containing the parameters of the planar figures to be displayed. We also provide a gallery of default phantoms of each type, accessible by calling the function \verb|my_phantomgallery|. 
The class \verb|Phantom| has two methods: \verb|get_phantom| that returns the image of the phantom, and \verb|get_sinogram| generates the analytical Radon transform, using a different procedure according to the phantom type. The implementation of the transforms is provided by the Scikit-image library \cite{scikitimage}.

The following lines of code provide an usage example of the \texttt{exact-sinogram} package.

\begin{lstlisting}[language=Python]
import numpy as np
from skimage.transform import radon, iradon

import exact_sinogram as es

circle = True
phantom_type = 'ellipses'  
#  phantom_type : 'ellipses' (or 'shepp_logan') , 'modified_shepp_logan', 'squares', 'rectangles'
n_points = 300;   # number of pixels

# Creation of an instance of the Phantom Class
Phm = es.Phantom(phantom_type = phantom_type, circle = circle)    

# Creation of the matrix-image of the phantom, w/ input the number of pixel:
P = Phm.get_phantom(N = n_points)

# Array of projection angles
theta_vec_deg = np.linspace(0, 359, 360)

# Exact sinogram and inversion (exact_sinogram package)
## Calculate the excact Sinogram
analytical_sinogram = Phm.get_sinogram(N = n_points, theta_vec = np.deg2rad(theta_vec_deg))

## invert the Radon transform on the exact sinogram
P_an = iradon(analytical_sinogram, theta = theta_vec_deg, circle=circle)

# Approximated sinogram and inversion (scikit-image)
## Calculate the approximated sinogram
approx_sinogram = radon(analytical_sinogram, theta = theta_vec_deg, circle=circle)
## invert the Radon transform on the approximated sinogram
P_approx = iradon(approx_sinogram, theta = theta_vec_deg, circle=circle)
\end{lstlisting}

\section{Results and conclusions}
For each phantom type, we compute the error between the synthetic data and the reconstructed data obtained by the following procedures:
\begin{enumerate}
    \item by computing the analytic sinogram and then applying the inverse Radon transform;
    \item in the classic way, by using the Radon transform and its inverse. 
\end{enumerate}

We are going to evaluate the quality of the reconstructed image by measuring the classical $2$-norm and we collect the results in Table \ref{tab:errors}.

\begin{table}[h!]
    \centering
    \begin{tabular}{|c|c|c|c|c|}
    \hline
    & Ellipses (Shepp-Logan) & Ellipses (Modified Shepp-Logan) & Squares & Rectangles \\ \hline
        $\frac{\Vert P_{an} - P \Vert_2}{\Vert P \Vert_2}$  &  0.03073 & 0.02590 & 0.02703 & 0.01838  \\ \hline
        $\frac{\Vert P_{inv} - P \Vert_2}{\Vert P \Vert_2}$  &  0.03604 & 0.03191 & 0.01997 & 0.02588 \\        
        \hline
    \end{tabular}
    \caption{Relative errors of the image reconstructed from the analytical sinogram $P_{an}$ and the image reconstructed form the numerical sinogram $P_{inv}$ with respect to the original image $P$, taking into account Gibbs phenomenon.}
    \label{tab:errors}
\end{table}

The Gibbs phenomenon \cite{Jerri_Gibbs} is the concentration of the errors along the discontinuities of the figure. With a mask, we eliminated the influence of boundaries and appreciated our algorithm that gives a smaller relative error (see Table \ref{tab:errors}) with respect to the reconstruction through the iradon function.

\bibliography{references}

\end{document}